\documentclass[10pt]{article}
\usepackage{amssymb,amsmath}
\usepackage[all]{xy}
\newtheorem{theorem}{Theorem}
\newtheorem{definition}{Definition}
\newtheorem{lemma}{Lemma}

\newtheorem{corollary}{Corollary}
\newtheorem{remark}{Remark}

\def\EE{{\mathbb E}}

\def\PP{{\mathbb P}}

\def\ind{{\mathbf{1}}}
\newcommand\bp[1]{\noindent {\em Proof{#1}.} $\quad$}
\def\ep{\hfill $\Box$}
\def\supp{{\mathrm{supp\,}}}

\def\m{{\bx}}
\def\mm{{\by}}
\def\cI{\Gamma}

\def\A{{\Gamma}}
\def\R{{\mathbb R}}

\def\N{{\mathbb N}}

\def\cG{{\mathcal G}}

\def\cF{{\mathcal F}}

\def\cP{{\mathcal P}}
\def\cB{{\mathcal B}}
\def\cL{{\mathcal L}}

\def\bw{{\mathbf w}}

\def\bx{\mathbf{x}}
\def\bxo{\overline{\mathbf{x}}}
\def\xo{\overline{x}}

\def\by{\mathbf{y}}

\title{The cavity method for counting spanning subgraphs subject to local constraints}
\author{Justin Salez\footnote{INRIA-\'Ecole
    Normale Sup\'erieure - France. Email: justin.salez@ens.fr}}
\begin{document}
\maketitle

\begin{abstract}
Using the theory of negative association for measures and the notion of random weak limits of sparse graphs, we establish the validity of the cavity method for counting spanning subgraphs subject to local constraints in asymptotically tree-like graphs. Specifically, the corresponding free entropy density is shown to converge along any sequence of graphs whose random weak limit is a tree, and the limit is directly expressed in terms of the unique solution to a limiting cavity equation. On a Galton-Watson tree, the latter simplifies into a recursive distributional equation which can be solved explicitely. As an illustration, we provide an explicit-limit formula for the $b-$matching number of an Erd\H{o}s-R\'enyi random graph with fixed average degree and diverging size, for any $b\in\N$.
\end{abstract}

\section{Introduction}

The general framework we consider is that of a finite graph $G=(V,E)$, in which \textit{spanning subgraphs} are weighted according to their local aspect around each vertex as follows : 
\begin{eqnarray}
\label{eq:productform}
\mu(F) & = & \prod_{i\in V}\mu_i(F\cap E_i).
\end{eqnarray}
Here, a spanning subgraph $(V,F)$ is identified with its egde-set $F\subseteq E$, and each $\mu_i$ is a given non-negative function over the subsets of $E_i:=\{e\in E; e\textrm{ is incident to }i\}$. We call $\mu$ the \textit{global measure} induced by the \textit{local measures} $\mu_i,{i\in V}$. Of particular interest in combinatorial optimization is the number 
\begin{eqnarray}
\label{eq:rankdef}
M(G) & = & \max\left\{|F|: F\in\supp(\mu)\right\},
\end{eqnarray}
which is the maximum possible size of a spanning subgraph $F$ satisfying the local contraint $\mu_i(F\cap E_i)>0$ at every node $i\in V$. More generally, counting the weighted number of spanning subgraphs of each given size in $G$, i.e. determining the generating polynomial
\begin{eqnarray}
\label{eq:pf}
Z(G;t) & = & \sum_{F\subseteq E}\mu(F)t^{|F|}
\end{eqnarray}
is a fundamental task, of which many combinatorial problems are special instances. Intimately related to this is the study of a random spanning subgraph $\cF$ sampled from the \textit{Gibbs-Boltzmann} law :
\begin{eqnarray}
\label{eq:boltzmann}
\PP^{t}_G(\cF=F) & = &\frac{\mu(F)t^{|F|}}{Z(G;t)},
\end{eqnarray}
where $t>0$ is a variable parameter called the \textit{activity}. In particular, the expected size of $\cF$ is called the \textit{energy} $U(G;t)$ and is connected to $Z(G;t)$ via the elementary identity
\begin{eqnarray}
\label{eq:logderivative}
U(G;t) & = & t\frac{d}{dt}\log Z(G;t).
\end{eqnarray}
Our concern is the behavior of these quantities in the \textit{infinite volume limit} : $|V|\to\infty$, $|E|=\Theta(|V|)$. 

Originating from spin glass theory \cite{mezardparisi}, the \textit{cavity method} is a powerful nonrigorous technique for evaluating such asymptotics on graphs that are \textit{locally tree-like}.
%, which is the case in the \textit{infinite volume limit} ($|V|\to\infty$) for many classical ensembles of diluted graphs. 
Essentially, the heuristic consists in neglecting cycles in order to obtain an approximate local fixed point equation for the marginals of the Gibbs-Boltzmann law. 
%$\mu^{(\lambda)}_G(\cF=F)$ amounts to c$\cF$, which would produce the correct answer if $G$ were indeed a tree. 
Despite its remarkable practical efficiency and the mathematical confirmation of its  analytical predictions for various important models \cite{talagrand,zeta,mwis,urbanke,bayatinair,dembomontanari,matcharxiv}, this \textit{ansatz} is still far from being completely understood, and the exact conditions for its validity remain unknown. More precisely, two crucial questions arise in presence of cycles :
\begin{enumerate}
\item \textbf{convergence} : is there a unique, globally attractive fixed
  point to the cavity equation ?
\item \textbf{correctness} : if yes, does it have any relation to the Gibbs-Boltzmann marginals ?
\end{enumerate}

In this paper, we exhibit a general condition under which the cavity method is valid for counting spanning subgraphs subject to local constaints. Specifically, we positively answer question 1 for arbitrary finite graphs (Theorem \ref{th:convergence}), under the only assumption that each local measure enjoys a certain form of negative association which we call the \textit{cavity-monotone} property, and which simply boils down to \textit{ultra-log-concavity} in the exchangeable case. Regarding question 2, we use the framework of \textit{local weak convergence} \cite{bensch,aldste} and the notion of \textit{unimodularity} \cite{aldlyo} to prove asymptotical correctness for any sequence of graphs whose random weak limit is concentrated on trees (Theorem \ref{th:main}). This includes many classical sequences of \textit{diluted graphs}, such as random $d-$regular graphs, Erd\H{o}s-R\'enyi random graphs with fixed average degree, or more generally random graphs with a prescribed degree distribution. 
In all these examples, the limit is a \textit{unimodular Galton-Watson (UGW) tree}. Thanks to the distributional self-similarity of such a tree, the cavity equation simplifies into a \textit{recursive distributional equation} which may be solved explicitely. As a motivation, let us first describe the implications of our work in the special case of \textit{b-matchings}. 

%\subsection{Applications to b-matchings}
A famous combinatorial structure that fits in the above framework is obtained by fixing $b\in\N$ and taking $\mu_i(F)=\ind(|F|\leq b)$ for all $i\in V$ : the induced global measure $\mu$ is then nothing but the counting measure for $b-$matchings in $G$, i.e. spanning subgraphs with maximum degree at most $b$. The reader is refered to the monograph \cite{schrijver} for a complete survey on $b-$matchings, and to \cite{lovaszplummer} for the important case of \textit{matchings} ($b=1$). The associated quantities $M_b(G)$ and $Z_b(G;t)$ are important graph invariants respectively known as the \textit{$b-$matching number} and \textit{$b-$matching polynomial}. Determining $Z_1(G;t)$ is a classical example of a computationally hard problem \cite{valiant}, although efficient approximation algorithms have been designed  \cite{stoc07,bangam}. The mathematical properties of $Z_b(G;t)$ have been investigated in detail, notably in the case $b=1$ for the purpose of understanding \textit{monomer-dimer systems} \cite{heilmannlieb,berg}. Interestingly, the geometry of the complex zeros of $Z_b(G,t)$ has been proven to be quite remarkable (see \cite{heilmannlieb} for $b=1$, \cite{ruelle} for $b=2$, and \cite{wagner} for the general case). Regarding $M_1(G)$, the first results in the infinite volume limit were obtained by Karp and Sipser \cite{karpsipser} for the Erd\H{o}s-R\'enyi random graph $G_n$ with average degree $c>0$ on $n$ vertices :
\begin{equation}
\label{eq:KS}
\frac{1}{n}M_1(G_n)\xrightarrow[n\to\infty]{P}1-\frac{t_c+e^{-ct_c}+ct_ce^{-ct_c}}{2},
\end{equation}
where $t_c\in(0,1)$ is the smallest root of $t=e^{-ce^{-ct}}$. The analysis has then been extended to random graphs with a log-concave degree profile \cite{friezebohman}, and finally to any graph sequence that converges in the local weak sense \cite{matcharxiv}. Contrastingly, only little is known for $b\geq 2$ : to the best of our knowledge, the limit of $\frac{1}{|V_n|}M_b(G_n)$ is only known to exist in the Erd\H{o}s-R\'enyi case \cite{gamarnik}, and could not be explicitely determined. As a special case of our main result, it will follow that
\begin{theorem}[$b-$matchings in locally tree-like graphs]
\label{th:KSb}
For any sequence of finite graphs $(G_n)_{n\in\N}$ satisfying $|E_n|=O(|V_n|)$ and whose random weak limit $\cL$ is concentrated on trees, the limits
\begin{eqnarray*}
f_b(\cL):=\lim_{n\to\infty}\frac{1}{|V_n|}\log Z_b(G_{n};t) & \textrm{ and } & m_b(\cL):=\lim_{n\to\infty}\frac{M_b(G_{n})}{|V_n|}
\end{eqnarray*}
exist and depend only on the random weak limit $\cL$. When $\cL$ is a UGW tree, we have the explicit formula
\begin{equation*}
 m_b(\cL):=\min_{s\in[0,1]}\left\{b-\frac{b}{2}g_b(s)-\frac{b}{2}(g_b\circ f_b)(s)+\frac{c}{2}f_b(s)(f_b\circ f_b)(s)\right\},
\end{equation*}
where $c,f,g$ are defined in terms of the degree generating function $\phi(s)=\sum_{k}\pi_ks^k$ as follows :
\begin{equation*}
c=\phi'(1),\qquad \textrm{  } \qquad f_b(s)=\frac{1}{c}\sum_{k=0}^{b-1}\frac{s^{k}\phi^{(k+1)}(1-s)}{k!}\qquad\textrm{ and }\qquad g_b(s)=\sum_{k=0}^b\frac{s^{k}\phi^{(k)}(1-s)}{k!}.
\end{equation*}
Moreover, any $s$ where the above minimum is achieved must be a root of $s=(f_b\circ f_b)(s)$.
\end{theorem}
For example, in the case of Erd\H{o}s-R\'enyi random graphs with average degree $c>0$ on $n$ vertices, the random weak limit $\cL$ is a.s. the law of a UGW tree with Poisson(c) degree distribution, and hence,
\begin{equation*}
\frac{1}{n}M_b(G_n)\xrightarrow[n\to\infty]{a.s.}\min_{s\in[0,1]}\left\{b-\frac{b}{2}\varphi_{b+1}(s) - \frac{b}{2}(\varphi_{b+1}\circ\varphi_{b})(s) + \frac{c}{2}\varphi_{b}(s)(\varphi_{b}\circ \varphi_{b})(s)\right\},
\end{equation*}
where we have set $$\varphi_b(s)=e^{-ct}\sum_{k=0}^{b-1}\frac{{cs}^k}{k!}.$$ Since any $s$ where the minimum is achieved must satisfy $s=(\varphi_{b}\circ \varphi_{b})(s)$, we recover exactly (\ref{eq:KS}) in the special case of matchings ($b=1$).

The paper is organized as follows : in section \ref{sec:prelim}, we recall the basic notions and properties pertaining to measures over subsets, which will be of constant use throughout the paper. In section \ref{sec:finite}, we define and study the cavity equation associated to a finite network. In section \ref{sec:infinite}, we extend the results to infinite networks that arise as local weak limits of finite networks. Finally, section  \ref{sec:activity} is devoted to the study of the cavity equation in the limit of infinite activity, and to its explicit resolution in the case of $b-$matchings.

\section{Preliminaries}
\label{sec:prelim}
In this section, we define the important notions pertaining to (non-negative) measures $\mu$ over the subsets of an arbitrary finite \textit{ground set} $E$. Later on, these will be specialized to the local measures $(\mu_i)_{i\in V}$ attached to the vertices of a graph $G$. First, $\mu$ is caracterized by its \textit{multivariate generating polynomial}
\begin{eqnarray}
\label{eq:gp}
Z(\bw) & = & \sum_{F\subseteq E}\mu(F)\bw^{F},
\end{eqnarray}
where $\bw=(w_e)_{e\in E}$ and $\bw^F=\prod_{e\in F}w_e$. Since $Z$ is affine in each $w_e,e\in E$, it can be decomposed  as 
\begin{eqnarray}
\label{eq:decomposition}
Z(\bw) & = & w_eZ^{/e}(\bw')+Z^{\setminus e}(\bw'),
\end{eqnarray}
where $\bw'=(w_f)_{f\neq e}$ and $Z^{\setminus e}$,  $Z^{/e}$ are the multi-affine polynomials with
ground set $E\setminus e$ respectively obtained from $Z$ by setting the
variable $w_e$ to $0$ (\textit{deletion}) and differentiating with respect to $w_e$ (\textit{contraction}). By definition, the \textit{cavity ratio} of the pair $(\mu,e)$ is then simply the multi-affine rational function
\begin{eqnarray}
\label{eq:cavity}
\cI^e_\mu(\bw') & = & \frac{Z^{/e}(\bw')}{Z^{\setminus e}(\bw')}.
\end{eqnarray}
When positive values are assigned to the variables (a so-called \textit{external field}), we may consider the probability distribution 
%\begin{eqnarray}
%\label{eq:boltzmann}
$\PP^\bw_\mu(\cF=F)  = {\mu(F)\bw^F}/{Z(\bw)}.$
%\end{eqnarray}
A quantity of interest is the expected size of $\cF$ when viewed as a function of the external field. We call this the \textit{energy} : 
\begin{eqnarray}
U_\mu(\bw) & = & \EE^\bw_\mu\left[|\cF|\right].
\end{eqnarray} 
From the decomposition (\ref{eq:decomposition}), it follows immediately that 
\begin{eqnarray}
\label{eq:energycavity}
\PP^\bw_\mu(e\in\cF) \ = \ \frac{w_e\cI^e_\mu(\bw')}{1+w_e\cI^e_\mu(\bw')}, & \textrm{ hence } & U_\mu(\bw) \ = \ \sum_{e\in E}\frac{w_e\cI^e_\mu(\bw')}{1+w_e\cI^e_\mu(\bw')}.
\end{eqnarray}
Note that the supremum of the energy is exactly the \textit{rank} of $\mu$ : $\rm{rank}(\mu)  =  \max\left\{|F|; F\in\supp(\mu)\right\}.$
The following properties will be of crucial importance throughout the paper.
\begin{definition}[Cavity-monotone measures]
The measure $\mu$ is called
\begin{itemize}
%\item $\mu$ is called \textit{deletion-tolerant} if its support is closed under inclusion : 
%$$\forall F'\subseteq F\subseteq E, \mu(F)>0\Longrightarrow\mu(F')>0.$$
\item \textit{Rayleigh} if every two distinct ground elements $e\neq f$ are negatively correlated in $\cF$ :
$$\forall \bw>0, \PP^\bw_\mu\left(e\in \cF,f\in \cF\right)\leq \PP^\bw_\mu\left(e\in \cF\right)\PP^\bw_\mu\left(f\in \cF\right).$$
\item \textit{Size-increasing} if every ground element $e$ positively influences the total size $|\cF|$ : 
$$\forall \bw>0, \EE^\bw_\mu\left[|\cF| \ind_{(e\in \cF)}\right]>\EE^\bw_\mu\left[|\cF|\right]\PP^\bw_\mu(e\in\cF).$$ 
\item \textit{Cavity-monotone} if its satisfies $\mu(\emptyset)>0$ and the two above properties. 
\end{itemize}
\end{definition}
Rayleigh measures were introduced in the context of matroid theory \cite{wagner08}, but soon found their place in the modern theory of negative dependence for probability measures \cite{pemantle, kahnneiman}. Cavity-monotone measures will play a major role in our study, for the following elementary reason. 
\begin{lemma}[Monotony of energy and cavity ratios]
\label{lm:equiv}
\begin{eqnarray*}
\mu(\emptyset)>0 & \Longleftrightarrow  & \textrm{ the cavity ratios }\cI^e_\mu,e\in E\textrm{ are well-defined on } [0,\infty)\times\ldots\times [0,\infty). \frac{}{}\\
\mu\textrm{ is Rayleigh } & \Longleftrightarrow & \textrm{ the cavity ratios } \cI^e_\mu,e\in E \textrm{ are non-increasing in each variable}.\frac{}{}\\
\mu\textrm{ is size-increasing} & \Longleftrightarrow & \textrm{ the energy } U_\mu \textrm{ is increasing in each variable} ;\\
& \Longleftrightarrow & \textrm{ for each } e\in E\textrm{ and }\bw>0, t\mapsto t\cI^e_\mu(t\bw')\textrm{ is increasing}.
\end{eqnarray*}
\end{lemma}
\bp{}
Differentiating the corresponding quantities and playing with the definition of $\PP_\mu^\bw$ easily yields
\begin{eqnarray*}
\frac{\partial \cI^e_\mu(\bw')}{\partial w_f} &  =  &
%\frac{Z^{/e/f}Z-Z^{/e}Z^{f/}}{Z^{\setminus e}Z^{\setminus e}}\ = \ 
%\frac{
\left(
\PP^\bw_\mu\left(e\in \cF,f\in \cF\right)- \PP^\bw_\mu\left(e\in \cF\right)\PP^\bw_\mu\left(f\in \cF\right)
\right)/
%}{
w_ew_f\PP^\bw_\mu\left(e\notin \cF\right)^2.\\
\frac{\partial U_\mu(\bw)}{\partial w_e} &  = &
%\frac{Z\sum_{F\ni e}|F|\mu(F)\bw^{F\setminus e}-Z^{/e}\sum_{F}|F|\mu(F)\bw^F}{Z^2}  \ = \ 
%\frac{
\left(\EE^\bw_\mu\left[|\cF| \ind_{(e\in \cF)}\right]-\EE^\bw_\mu\left[|\cF|\right]\PP^\bw_\mu(e\in\cF)
\right)/
%}{
w_e.\\
\frac{\partial t\Gamma^e_\mu(t\bw')}{\partial t} &  = & 
%\frac{
\left(
\EE^\bw_\mu\left[|\cF| \ind_{(e\in \cF)}\right]-\EE^\bw_\mu\left[|\cF|\right]\PP^\bw_\mu(e\in\cF)
\right)/
%}{
tw_e\PP^\bw_\mu\left(e\notin \cF\right)^2. 
\end{eqnarray*}

\begin{remark}[Matroids]
\label{re:support}
Interestingly, the support of a cavity-monotone measure admits a remarkable structure : it follows from \cite[Theorem 4.6]{wagner08} that for $\mu$ Rayleigh with $\mu(\emptyset)>0$, $I=\supp \mu$ is a \textit{matroid}:
\begin{itemize}
\item $I$ is not empty ;
\item If $B\in I$ and $A\subseteq B$, then $A\in I$ ;
\item If $A,B\in I$ and $|A| < |B|$, then $\exists e\in B\setminus A$ such that $A\cup e\in I$.
\end{itemize}
\end{remark}

The cavity-monotone property admits a particularly simple caracterization in the important case where $\mu$ is \textit{exchangeable}, i.e. $\mu(F)=c(|F|)$ for some non-negative coefficients $c(0),\ldots,c(m), m=|E|$:
\begin{lemma}[The exchangeable case]
\label{lm:exchangeable}
An exchangeable measure $\mu$ is cavity-monotone if and only if 
\begin{enumerate}
\item $c$ is log-concave, i.e. $c^2(k)\geq c(k-1)c(k+1)$ for all $0<k<m$, and 
\item the support $\{0\leq k\leq m:c(k)>0\}$ is an interval containing $0$ and $1$.
\end{enumerate}
\end{lemma}
In particular, so is the measure $\mu(F)=\ind_{(|F|\leq b)}$ describing the local constraints of a $b-$matching.

\bp{ of Lemma \ref{lm:exchangeable}}
The result essentially follows from the work of Pemantle \cite{pemantle}. Indeed, Theorem 2.7 therein guarantees that $\mu$ is Rayleigh if and only if the sequence $c$ is log-concave and its support is an interval. That the latter must contain $0$ is nothing but the last property in the definition of a cavity-monotone measure. That it is not reduced to $0$ is imposed by the strict inequality in the size-increasing property. Conversely, let us show that any exchangeable measure $\mu$ with $c(0)>0$ and $c(1)>0$ is indeed size-increasing. Fix an external field $\bw>0$. By Lemma 2.9 in \cite{pemantle}, the law obtained from $\PP^\bw_\mu$ by conditionning on the event $\{|\cF|=k\}$ is stochastically increasing in $k$. By Proposition 1.2 in \cite{pemantle}, this implies in particular that for every $e\in E$, the following weak inequality holds :
$$\EE^\bw_\mu\left[|\cF| \big| e\in \cF\right]\geq \EE^\bw_\mu\left[|\cF|\right].$$
Note that the condition $c(1)>0$ guarantees that this conditional expectation is well-defined. Since we have not yet used the fact that $c(0)>0$, the above inequality remains true if one changes the coefficient $c(0)$ to $0$. Setting it then back to its initial (positive) value does not affect the left-hand side, but strictly decreases the right-hand side, hence the desired strict inequality. 
\ep

\section{The cavity equation on finite networks}
\label{sec:finite}
Let $G=(V,E)$ be a finite graph at the vertices of which some local measures $\mu_i,i\in V$ are specified. We call the resulting object a \textit{network}. A \textit{configuration} $\m$ is an assignment of numbers $x_{i\to j}\geq 0$ to every oriented edge $i\to j\in\vec{E}$. Starting from a configuration $\m$, we define a new configuration $\mm=\A_G(\m)$ by 
\begin{equation}
\label{eq:BPdef}
y_{i\to j} = \cI^{ij}_{\mu_i}\left(x_{k\to i} : k\in\partial i\setminus j\right),
\end{equation}
where $\partial i$ denotes the set of all neighbors of $i$. Each $x_{i\to j}$ may be thought of as a message sent by $i$ to $j$ along the edge $ij$, and $\A_G$ as a local rule for propagating messages. For $t>0$, the fixed point equation 
\begin{equation}
\label{eq:cavityequation}
\m=t\A_G(\m)
\end{equation} is called the \textit{cavity equation} at activity $t$ on the network $G$. Its relation to the global measure $\mu$ induced by the  $(\mu_i)_{i\in V}$ is revealed by the following well-known result.
\begin{lemma}[Validity on trees]
\label{lm:trees}
Assume that $G$ is finite and acyclic. Then, for every activity $t>0$,
\begin{enumerate}
\item \textbf{convergence} : the cavity equation admits a unique solution $\m(t)$, which can be reached from any initial configuration by iterating $t\A_G$ a number of times equal to the diameter of $G$ ;
\item\textbf{correctness} : for every $i\in V$, the exact marginal law of $\cF\cap E_i$ under the Gibbs-Boltzmann law $\PP^{t}_G$ is given by directly  imposing the external field $\{x_{j\to i}(t):j\in\partial i\}$ onto the local measure $\mu_i$. 
\end{enumerate}
\end{lemma}
The important consequence is that on trees, the energy $U(G;t)$ can be determined using only local operations : 
\begin{eqnarray}
\label{eq:cavityformula}
U(G;t) & = & \frac{1}{2}\sum_{i\in V}U_{\mu_i}(x_{j\to i}(t) : j\in\partial i) \ =  \ \sum_{ij\in E}\frac{x_{j\to i}(t)x_{i\to j}(t)}{t+x_{j\to i}(t)x_{i\to j}(t)},
\end{eqnarray}
where the second equality is obtained by applying (\ref{eq:energycavity}) to each $\mu_i,i\in V$. 

\bp{ of Lemma \ref{lm:trees} } When $i$ is a leaf, the message $y_{i\to j}$ defined by equation (\ref{eq:BPdef}) does not depend at all on the initial configuration $\bx$. Iterating this argument immediately proves the convergence part, and we now focus on correctness. Let $G=(V,E)$ be a finite tree, $\circ$ a vertex, and $i$ a neighbour of $\circ$. We let $G_{i\to\circ}$ denote the subtree induced by $\circ$ and all vertices that the edge $i\circ$ separates from $\circ$. Now assume that $G$ is equipped with local measures, and let $G_{i\to\circ}$ inherit from these local measures, except for $\mu_\circ$ which we replace by the trivial local measure with constant value $1$. With these notations, any spanning subgraph $F\subseteq E$ can be uniquely decomposed as the disjoint union of a subset $I\subseteq E_\circ$ and a spanning subgraph $F_i$ on each $G_{i\to\circ},i\in\partial\circ$,  with $i\circ\notin F_i$. Thus, writing $\mu_G$ for the global measure on the network $G$, we have for any $t>0$,
\begin{eqnarray*}
\mu_G(F)t^{|F|} & = & \mu_\circ(I)\prod_{i\in I}t^{|F_i|+1}\mu_{G_{i\to\circ}}({F_i\cup {i\circ}})\prod_{i\notin I}t^{|F_i|}\mu_{G_{i\to\circ}}({F_i}).\\
\end{eqnarray*}
Fixing $I$ and summing over all possible values for $F_i,i\in\partial\circ$, we obtain
\begin{eqnarray*}
\PP^t_G(\cF\cap E_\circ = I) & = & C \mu_\circ(I)\prod_{i\in I}\PP^t_{G_{i\to\circ}}({i\circ\in\cF})\prod_{i\notin I}\PP^t_{G_{i\to\circ}}({i\circ\notin\cF}).\\
& = & C'\mu_\circ(I)\prod_{i\in I}\frac{\PP^t_{G_{i\to\circ}}({i\circ\in\cF})}{\PP^t_{G_{i\to\circ}}({i\circ\notin\cF})},
\end{eqnarray*}
where $C,C'$ are normalizing constants that do not depend on $I$. This already proves that the law of $\cF\cap E_\circ$ can be  obtained from the local measure $\mu_\circ$ by imposing on each edge $i\circ\in E_\circ$ the external field 
\begin{equation}
\label{eq:treerec}
x_{i\to\circ}(t):=\frac{\PP^t_{G_{i\to\circ}}({i\circ\in\cF})}{\PP^t_{G_{i\to\circ}}({i\circ\notin\cF})}.
\end{equation}
In turn, this ratio can now be computed by applying the result to the vertex $i$ in the network $G_{i\to\circ}$ :
\begin{eqnarray*}
\frac{\PP^t_{G_{i\to\circ}}(i\circ \in \cF)}{\PP^t_{G_{i\to\circ}}(i\circ \notin \cF)} & =  & t\Gamma^{i\circ}_{\mu_i}\left(x_{k\to i}(t): k\in\partial i\setminus \circ\right),
\end{eqnarray*}
which shows that the configuration $\bx(t)$ defined on $G$ by (\ref{eq:treerec}) satisfies the cavity equation (\ref{eq:cavityequation}).
\ep

There are two distinct parts in Lemma \ref{lm:trees} : convergence and correctness. As we will now show, the former extends to arbitrary graphs under the only  assumption that each local measure is cavity-monotone. Henceforth, such a network  will be called a \textit{cavity-monotone network}. 

\begin{theorem}[Convergence on finite cavity-monotone networks]
\label{th:convergence}
On a finite cavity-monotone network, the cavity equation admits a unique, globally attractive fixed point $\bx(t)$ at any activity $t>0$. 
\end{theorem}

\bp{}
Fixing $t>0$ and starting with the minimal configuration $\m^0:=\bf{0}$, we set inductively 
$$\m^{k+1}(t):=t\A_G(\m^k(t)),$$
for all $k\in\N$. By Lemma \ref{lm:equiv}, the Rayleigh property of the local measures $\mu_i,i\in V$ ensures that $\A_G$ is coordinate-wise non-increasing on the space of configurations. Therefore, the limiting configuration 
\begin{equation}
\m^{-}(t):=\lim_{k\to\infty}\uparrow\m^{2k}(t)\qquad\textrm{ and }\qquad \m^{+}(t):=\lim_{k\to\infty}\downarrow\m^{2k+1}(t)
\end{equation}
exist, and any fixed point $\m=t\A_G(\m)$ must satisfy $\m^{-}(t)\leq\m\leq\m^+(t).$ Moreover, $\A_G$ is clearly continuous with respect to the product topology on configurations, so that $t\A_G(\m^-(t))=\m^+(t)$ and $t\A_G(\m^+(t))=\m^-(t)$. Thus, the existence of unique globally attractive solution to (\ref{eq:cavityequation}) boils down to the equality 
\begin{equation}
\label{eq:toshow}
\m^-(t)=\m^+(t). 
\end{equation}
Now applying (\ref{eq:energycavity}) to the local measure at a fixed vertex $i\in V$ yields 
\begin{eqnarray*}
U_{\mu_i}\left(x^-_{j\to i}(t): j\in\partial i\right) \ = \ \sum_{j\in\partial i}\frac{x^-_{j\to i}(t)x^+_{i\to j}(t)}{t+x^-_{j\to i}(t)x^+_{i\to j}(t)},
&  & U_{\mu_i}\left(x^+_{j\to i}(t): j\in\partial i\right) \ = \ \sum_{j\in\partial i}\frac{x^+_{j\to i}(t)x^-_{i\to j}(t)}{t+x^+_{j\to i}(t)x^-_{i\to j}(t)}.
\end{eqnarray*}
Summing over all vertices $i\in V$, we therefore obtain 
\begin{eqnarray*}
\sum_{i\in V}U_{\mu_i}\left(x^-_{j\to i}(t): j\in\partial i\right)&  = & 
\sum_{ij\in E}\left(\frac{x^-_{j\to i}(t)x^+_{i\to j}(t)}{t+x^-_{j\to i}(t)x^+_{i\to j}(t)}+\frac{x^+_{j\to i}x^-_{i\to j}(t)}{t+x^+_{j\to i}(t)x^-_{i\to j}(t)}\right)\\
& = & \sum_{i\in V}U_{\mu_i}\left(x^+_{j\to i}(t): j\in\partial i\right).
\end{eqnarray*}
This implies (\ref{eq:toshow}), since by Lemma \ref{lm:equiv} each $U_{\mu_i},i\in V$ is strictly increasing in every coordinate.
\ep

\section{The limit of infinite volume}
\label{sec:infinite}
In the previous section, we have established existence and uniqueness of a \textit{cavity solution} on any finite cavity-monotone network. Our concern now is its asymptotical meaning as the size of the underlying graph tends to infinity. Following the principles of the \textit{objective method} \cite{aldste}, we will replace the asymptotical analysis of our finite networks by the direct study of their infinite limits. 
\subsection{Random weak limits}
We first briefly recall the framework of local convergence, introduced by Benjamini and Schramm  \cite{bensch} and developped further by Aldous and Steele \cite{aldste}. Examples of successful uses include \cite{aldlyo, BenjaminiSS08,elek,lyotree,borgs2010left,rankarxiv,matcharxiv,dembomontanari}. Here, a \textit{network} will be simply a denumerable graph $G=(V,E)$ whose vertices are equipped with local measures $\mu_i,i\in V$. A \textit{rooted network} $(G,\circ)$ is a network together with the specification of a particular vertex $\circ\in V$, called the
\textit{root}. For $\varepsilon\geq 0$, we write $(G',\circ')\stackrel{\varepsilon}{\equiv}(G,\circ)$ if there exists a bijection $\gamma\colon V\to V'$ that preserves
\begin{itemize}
\item the root : $\gamma(\circ)=\circ'$ ;
\item the adjacency : $ij\in E\Longleftrightarrow \gamma(i)\gamma(j)\in E'$ ;
\item the support of the local measures : $\mu_i(F)>0 \Longleftrightarrow \mu'_{\gamma(i)}(\gamma(F))>0$, with $\gamma(F)=\{\gamma(i)\gamma(j):ij\in F\}$.
\item the values of the local measures, up to $\varepsilon$ : $|\mu'_{\gamma(i)}(\gamma(F))-\mu_i(F)|\leq \varepsilon$.
\end{itemize}
We let $\cG_*$ denote the set of all locally finite connected rooted networks considered up to the isomorphism relation $\stackrel{0}{\equiv}$. In the space $\cG_*$, a sequence $\left\{(G_n,\circ_n);n\in \N\right\}$ \textit{converges locally} to $(G,\circ)$ if for every radius $k\in\N$ and every $\varepsilon>0$, there is $n_{k,\varepsilon}\in\N$ such that
$$n\geq n_k\Longrightarrow [G_n,\circ_n]_k\stackrel{\varepsilon}{\equiv}[G,\circ]_k,$$ 
where $[G,\circ]_k$ denotes the finite rooted network obtained by keeping only the vertices lying at graph-distance at most $k$ from $\circ$. 
It is not hard to construct a distance which metrizes this notion of convergence and turns $\cG$ into a complete separable metric space. We can thus import the usual machinery of weak convergence of probability measures on Polish spaces (see e.g. \cite{billingsley}).
%Rather than just graphs $G=(V,E)$, it will be sometimes convenient to
%work with \textit{discrete networks} $G=(V,E,\cM)$, in which the additional
%specification of a \textit{mark map} $\cM\colon E\to\N$ allows to attach useful local
%information to edges, such as their absence/presence in a certain matching. We then simply require the isomorphisms in the
%above definition to preserve these marks.

\textit{Uniform rooting} is a natural procedure for turning a finite deterministic
network $G$ into a random element of $\cG_*$ : one simply chooses uniformly at random a vertex
$\circ$ to be the root, and restrains $G$ to the
connected component of $\circ$. If $(G_n)_{n\in\N}$ is a sequence of finite networks and if the sequence of their laws under uniform rooting admits a weak limit $\cL\in\cP(\cG_*)$, we call $\cL$ the \textit{random weak limit} of the sequence  $(G_n)_{n\in\N}$. In \cite{aldlyo}, it was shown that any such limit enjoys a remarkable invariance property known as \textit{unimodularity} : let $\cG_{**}$ denote the space of locally finite connected networks with an ordered pair of distinguished adjacent vertices $(G,\circ,i)$, taken up to the natural isomorphism relation and endowed with the natural topology. A measure $\cL\in\cP(\cG_{*})$ is called unimodular if it satisfies the \textit{Mass-Transport Principle} : for any Borel function $f\colon \cG_{**}\to [0,\infty]$, 
\begin{equation}
\label{eq:mtp}
\cL\left[\sum_{i\in\partial\circ}f(G,\circ,i)\right]=\cL\left[\sum_{i\in\partial\circ}f(G,i,\circ)\right],
\end{equation}
where we have written $\cL[\cdot]$ for the expectation with respect to $\cL$.
This is a deep and powerful notion, which we will now use to extend the results of section \ref{sec:finite} to the infinite setting. 

\subsection{Main result : validity of the cavity method on unimodular trees}
The definition of $\A_G$ remains valid for any locally finite network $G$. When the latter is cavity-monotone, the configurations $\m^{-}(t),\m^+(t)$ introduced in the proof of Theorem \ref{th:convergence} remain perfectly well-defined, and the convergence of the cavity method again boils down to the identity $\bx^-(t)=\bx^+(t)$. However, the proof of the latter involves a summation over all edges, which is no longer valid in the infinite setting. Instead, the desired $\bx^-(t)=\bx^+(t)$ will be derived from unimodularity, and will thus hold for any random weak limit of finite networks. Indeed, applying the Mass-Transport Principle to the function 
$$f(G,\circ,i) := \frac{x^-_{i\to\circ}(t)x^+_{\circ\to i}(t)}{t+x^-_{i\to\circ}(t)x^+_{\circ\to i}(t)}$$
yields $\cL\left[U_{\mu_\circ}(x^-_{i\to\circ}(t):i\in\partial\circ)\right] = \cL\left[U_{\mu_{\circ}}(x^+_{i\to\circ}(t):i\in\partial\circ)\right]$ ($f$ is Borel as the pointwise limit of continuous functions). Under the assumption $\cL\left[\rm{rank}\,(\mu_\circ)\right]<\infty$, this expectation is finite, and the size-increasing property of $\mu_\circ$ then implies that $\cL-$almost surely, $x^-_{i\to\circ}(t)=x^{+}_{i\to\circ}(t)$ for all $i\in\partial\circ$. This automatically extends to every oriented edge since under unimodularity, \textit{everything shows up at the root} (another fruitful application of the Mass-Transport-Principle, see \cite[Lemma 2.3]{aldlyo}). We state this as a Theorem.

\begin{theorem}[Convergence of the cavity method on unimodular networks]
\label{th:unimodular}
Let $\cL$ be a unimodular probability measure supported by cavity-monotone networks. If $\cL\left[\rm{rank}\,(\mu_\circ)\right]<\infty$, then the cavity equation admits $\cL-$a.-s. a unique, globally attractive solution $\bx(t)$ at any activity $t>0$. 
\end{theorem}

By analogy with formula (\ref{eq:cavityformula}) in the finite case, the (now well-defined) quantity
\begin{equation}
\label{eq:ulim}
u(\cL;t)=\frac{1}{2}\cL\left[\sum_{i\in\partial\circ}\frac{x_{i\to\circ}(t)x_{\circ\to i}(t)}{t+x_{i\to\circ}(t)x_{\circ\to i}(t)}\right]%U_{\mu_\circ}\left(x_{i\to \circ}(t): i\sim\circ\right)\right].
\end{equation}
appears as a natural candidate for the limiting energy of any sequence of finite networks whose random weak limit is $\cL$. Our second result is the validity of this \textit{cavity ansatz} when $\cL$ is concentrated on trees.

\begin{theorem}[Asymptotical correctness of the cavity method]
\label{th:main}
Let $(G_n)_{n\in\N}$ be a sequence of finite cavity-monotone networks admitting a random weak limit $\cL$ which is concentrated on cavity-monotone trees. Assume that the local rank at a uniformly chosen vertex is uniformly integrable as $n\to\infty$. Then, 
%\begin{equation*}
%\frac{\EE_{G_n}\left[|\cF|\right]}{|V_n|}\xrightarrow[n\to\infty]{}\frac{1}{2}\cL\left[\sum_{i\sim\circ}\frac{x_{i\to\circ}x_{\circ\to i}}{1+x_{i\to\circ}x_{\circ\to i}}\right],
%\end{equation*}
%where $\bx$ denotes the almost-surely unique solution to the cavity equation under $\cL$.
% \end{theorem}
%\begin{corollary}
%Under the abnove assumptions, we have 
\iffalse
\begin{eqnarray}
\label{eq:energy}
 \frac{U(G_n;t)}{|V_n|} & \xrightarrow[n\to\infty]{} & u(\cL;t) \\
\label{eq:ldp}
 \frac{1}{|V_n|}\log \frac{Z(G_n;t)}{Z(G_n;1)} & \xrightarrow[n\to\infty]{} & \int_{1}^t \frac{u(\cL;s)}s ds.
\end{eqnarray}
\fi
\begin{eqnarray}
\label{eq:energy}
 \frac{U(G_n;t)}{|V_n|} & \xrightarrow[n\to\infty]{} & u(\cL;t).
 \end{eqnarray}
If $|E_n|=O(|V_n|)$ and all the local measures take values in $\{0\}\cup K$ for a fixed compact $K\subseteq (0,\infty)$, then
\begin{eqnarray}
\label{eq:number}\frac{1}{|V_n|}\log {Z(G_n;t)} & \xrightarrow[n\to\infty]{} & \cL[\log\mu_\circ(\emptyset)] + \int_{0}^t \frac{u(\cL;s)}s ds,\\
\label{eq:rank}
\frac{M(G_n)}{|V_n|} & \xrightarrow[n\to\infty]{} & m(\cL):=\lim_{t\to\infty}\uparrow u(\cL;t).
\end{eqnarray}
\end{theorem}
%The important quantity $m(\cL)$ will be studied in details in the next section.
\begin{remark}[Large deviation principle]
\label{re:ldp}
Integrating (\ref{eq:energy}) immediately implies that 
\begin{eqnarray*}
\label{eq:ldp}
 \frac{1}{|V_n|}\log \frac{Z(G_n;t)}{Z(G_n;1)} & \xrightarrow[n\to\infty]{} & \int_{1}^t \frac{u(\cL;s)}s ds.
\end{eqnarray*}
It will later be checked that $t\mapsto u(\cL;t)$ is continuous on $\R_+$ (Remark \ref{re:continuity}). Therefore, denoting by $\cF_n$ a random spanning subgraph with the Gibbs-Boltzmann law $\PP^1_{G_n}$, G\"artner-Ellis Theorem \cite{dembozeitouni} guarantees that $|\cF_n|/|V_n|$ obeys a large deviation principle with rate $|V_n|$ and  good rate function $y\mapsto \int_{0}^\infty \left(y-u(\cL;e^s)\right)^+ ds$.
\end{remark}

\subsection{Proof of the main result}

\begin{lemma}[Tree approximation]
\label{lm:treeapprox}
Let $(G,\circ)$ be a finite rooted cavity-monotone network, and let $k\in\N$. If $[G,\circ]_{2k+2}$ is a tree, then for every activity $t>0$,
\begin{eqnarray*}
U_{\mu_\circ}\left(x^{2k}_{i\to\circ}(t):i\in\partial\circ\right)\leq \EE^t_G\left[|\cF\cap E_\circ|\right] \leq U_{\mu_\circ}\left(x^{2k+1}_{i\to\circ}(t):i\in\partial\circ\right).
\end{eqnarray*}
\end{lemma}
\bp{}
The proof makes use of a classical ingredient known as the \textit{spatial Markov property}, which we first briefly recall. Let $G=(V,E)$ be a finite network and let $S$ be an induced subgraph. We let $\partial S$ denote the boundary of $S$, i.e. the set of edges having one end-point in $S$ and one in $S^c$. Any \textit{boundary condition} $B\subseteq \partial S$ can be used to assign local measures to the vertices of $S$, namely $\mu^B_i(F):=\mu_i\left(F\cup (B\cap E_i)\right)$. Note that these local measures differ from the original ones only for vertices that are adjacent to the boundary. The resulting network is denoted by $S|B$. Now, a spanning subgraph $F\subseteq E$ is clearly the disjoint union of a spanning subgraph $F_{\rm{int}}$ of $S$, a boundary condition $B\subseteq \partial S$ and a spanning subgraph $F_{\rm{ext}}$ in $S^c$. The product form of $\mu_G$ immediately yields :
\begin{equation}
\PP^t_G(\cF=F)=\PP^t_{S|B}(\cF=F_{\rm{int}})\PP^t_G(\cF\cap\partial S=B)\PP^t_{S^c|B}(\cF=F_{\rm{ext}}).
\end{equation}
In other words, conditionally on the boundary $\cB:=\cF\cap \partial S$, the restrictions of $\cF$ to $S$ and $S^c$ are independent with law $\PP^t_{S|\cB}$ and $\PP^t_{S^c|\cB}$, respectively. Applying this to the tree $S=[G,\circ]_{2k+2}$,
\begin{eqnarray*}
\EE^t_G\left[|\cF\cap E_\circ|\right]& = &\sum_{B\subseteq \partial S}\PP^t_G(\cF\cap\partial S=B)\EE^t_{S|B}\left[|\cF\cap E_\circ|\right]\\
& = &  \sum_{B\subseteq \partial S}\PP^t_G(\cF\cap\partial S=B)U_{\mu_\circ}\left(x^{(B)}_{i\to\circ}(t):i\in\partial\circ\right),
\end{eqnarray*}
where we have applied Lemma \ref{lm:trees} to the tree $S|B$, writing $\bx^{(B)}(t)$ for the unique solution to the cavity equation at activity $t$ thereon. But by monotony of the cavity operator, each $x^{(B)}_{i\to\circ}(t),i\in\partial\circ$ must satisfy $x^{2k}_{i\to\circ}(t)\leq x^{(B)}_{i\to\circ}(t)\leq x^{2k+1}_{i\to\circ}(t)$. Using the size-increasing property of $\mu_\circ$, we see that 
$$U_{\mu_\circ}\left(x^{2k}_{i\to\circ}(t):i\in\partial\circ\right) \leq U_{\mu_\circ}\left(x^{(B)}_{i\to\circ}(t):i\in\partial\circ\right)\leq U_{\mu_\circ}\left(x^{2k+1}_{i\to\circ}(t):i\in\partial\circ\right),$$
and re-injecting this into the above equation finally yields the desired inequalities.
\ep

Let us now see how Lemma \ref{lm:treeapprox} implies the convergence (\ref{eq:energy}).  Let $(G_n)_{n\in\N}$ be a sequence of finite cavity-monotone networks admitting a random weak limit $\cL$ which is concentrated on cavity-monotone trees. Denote by $\cL_n\in\cP(\cG_*)$ the law under uniform rooting of $G_n$, so that $\cL_n\Rightarrow\cL$. We will use the short-hand $u_{k}(G,\circ)=U_{\mu_\circ}\left(x^{k}_{i\to\circ}(t):i\in\partial\circ\right)$, and $\chi_k(G,\circ)$ for the indicator function that $[G,\circ]_{2k+2}$ is a tree. Lemma \ref{lm:treeapprox} guarantees that for any finite cavity-monotone network $G$ and any vertex $\circ$,
\begin{equation}
\label{eq:ineqq}\chi_k(G,\circ)u_{2k}(G,\circ)\leq \EE^t_{G}\left[|\cF|\cap E_\circ\right]\leq \chi_k(G,\circ)u_{2k+1}(G,\circ) + (1-\chi_k(G,\circ))\rm{rank(\mu_\circ)}.
\end{equation}
As functions of $(G,\circ)$, the left-hand side and right-hand side are continuous on $\cG_*$, since they depend only on $[G,\circ]_{2k+2}$. Moreover, both are dominated by $(G,\circ)\mapsto \rm{rank(\mu_\circ)}$ which is assumed to be uniformly integrable with respect to the sequence $(\cL_n)_{n\in\N}$. Thus, their expectation under $\cL_n$ tends to their expectation under $\cL$ as $n\to\infty$. But $\chi_k$ is zero on the support of $\cL$, so we are simply left with 
\begin{eqnarray*}
\frac{1}{2}\cL\left[U_{\mu_\circ}\left(x^{2k}_{i\to\circ}(t):i\in\partial\circ\right)\right] 
\ \leq \ 
\liminf_{n\to\infty}\frac{U(G_n;t)}{|V_n|} 
& \leq & 
\limsup_{n\to\infty}\frac{U(G_n;t)}{|V_n|} 
\ \leq \ 
\frac{1}{2}\cL\left[U_{\mu_\circ}\left(x^{2k+1}_{i\to\circ}(t):i\in\partial\circ\right)\right].
\end{eqnarray*}
Since the random weak limit $\cL$ is unimodular, Theorem \ref{th:unimodular} finally implies that both the lower and upper bounds tend to $\cL\left[U_{\mu_\circ}\left(x_{i\to\circ}(t):i\in\partial\circ\right)\right]=u(\cL;t)$ as $k\to\infty$. Note that the  requirement $\cL[\rm{rank}(\mu_\circ)]<\infty$ in Theorem \ref{th:unimodular} is here automatically fullfilled, by the uniform integrability assumption. 

It now remains to show (\ref{eq:number}) and (\ref{eq:rank}). The identity (\ref{eq:logderivative}) implies that for any activity $t>0$ and any finite network $G$ satisfying $\mu(\emptyset)>0$, 
$$\frac{1}{|V|}\log {Z(G;t)} = \frac{1}{|V|}\sum_{\circ\in V}\log \mu_\circ(\emptyset) + \int_{0}^t \frac{U(G;s)}{s|V|} ds.$$ 
Now take $G=G_n$ and let $n\to\infty$ : the compactness assumption guarantees that $\log \mu_\circ(\emptyset)$ is bounded uniformly in $n$, so the first term converges to $\cL[\log\mu_\circ(\emptyset)]$. As per the second one, it tends to $\int_{0}^t \frac{u(\cL;s)}s ds$ because of (\ref{eq:energy}), provided the uniform domination holds in Lebesgue's dominated convergence Theorem. The latter fact is ensured by the first inequality in Lemma \ref{lm:unifctrl} below, combined with the compactness assumption and the fact that $|E_n|=O(|V_n|)$. The second inequality in Lemma \ref{lm:unifctrl} easily guarantees (\ref{eq:rank}).

\begin{lemma}[Uniform controls for the energy]
\label{lm:unifctrl}
Let $G$ be a finite cavity-monotone network. As a function of the activity $t$, the energy $U(G;t)$ increases from $0$ to $M(G)$. Furthermore, the rate of convergence to these two extrema can be precisely controlled :
\begin{eqnarray}
\label{eq:ctrlzero}
\forall t >0, \qquad U(G;t) & \leq & t \sum_{ij\in E}A(\mu_i)A(\mu_j);\\
\label{eq:ctrlinfinity}
\forall t > 1, \qquad U(G;t) & \geq & M(G) - \frac{1}{\log t}\left(|E|\log 2 + \sum_{i\in V}\log A(\mu_i)\right).
\end{eqnarray}
where $A(\mu)=\frac{\max \mu}{\min \mu}$, with $\max \mu={\max\{\mu(F):F\in\supp(\mu)\}}$ and $\min\mu ={\min\{\mu(F):F\in\supp(\mu)\}}$.
\end{lemma}

\bp{ of Lemma \ref{lm:unifctrl}}
That the energy increases with the activity is equivalent by (\ref{eq:logderivative}) to the convexity of $\theta\mapsto \log Z(G;e^\theta)$, a direct consequence of H\"older's inequality. This also implies that for any $t>1$
\begin{eqnarray*}
{U(G;t)}\log t & \geq & \log\frac{Z(G;t)}{Z(G;1)}.
\end{eqnarray*}
Clearly $Z(G;t)\geq t^{M(G)}\min \mu$ and $Z(G;1)\leq 2^{|E|}\max \mu$ and $A(\mu)\leq\prod_{i\in V}A(\mu_i)$,  so (\ref{eq:ctrlinfinity}) follows. Regarding (\ref{eq:ctrlzero}), we have for any $t>0$
\begin{eqnarray*}
{U(G;t)} & = & \sum_{ij\in E}{\PP^t(ij\in\cF)}\\
& = & t\sum_{ij\in E}\frac{\sum_{F\subseteq E\setminus ij}\mu(F\cup ij)t^{|F|}}{\sum_{F\subseteq E}\mu(F)t^{|F|}}\\
& \leq & t\sum_{ij\in E}\max_{F\subseteq E\setminus ij,\mu(F)>0}\frac{\mu(F\cup ij)}{\mu(F)}\\
& \leq & t\sum_{ij\in E}A(\mu_i)A(\mu_j),
\end{eqnarray*}
where the third line uses the standard inequality $\frac{a+b}{c+d}\leq \max(\frac{a}{c},\frac{b}{d})$ for any $a,b,c,d>0$, and the crucial fact that $\mu(F\cup ij)>0 \Longrightarrow \mu(F)>0$ (see Remark \ref{re:support}). 

\section{The limit of infinite activity}
\label{sec:activity}
The goal of this final section is to describe the important quantity $m(\cL)$ introduced in Theorem \ref{th:main} directly in terms of a certain local equation which we naturally call the \textit{cavity equation at infinite activity}. This will then be used to establish the explicit formulae that have been mentioned in the introduction.

\subsection{The cavity equation at infinite activity}

%In particular, we may consider the limiting cavity-ratio $\overline{\Gamma}^e_\mu(\bw):=\lim_{t\to\infty}\uparrow t\Gamma^e_\mu(t\bw),$ and the associated cavity-operator on any cavity-monotone network $G$ :
Let $G$ be a cavity-monotone network. From Lemma \ref{lm:equiv}, it follows that $(t,\bx)\mapsto t\Gamma_G(t\bx)$ is increasing in $t$ and decreasing in $\bx$. We may thus define a limiting cavity-operator by
$$\overline{\Gamma}_G(\bx):=\lim_{t\to\infty}\uparrow t\Gamma_G(t\bx).$$
By monotony, $\overline\Gamma_G\colon [0,\infty)^{\vec E}\to (0,\infty]^{\vec E}$ is well-defined without any ambiguity regarding the order in which the limits $t\to\infty$ and $\bx\to 0$ are taken. 
Note also that $\overline\Gamma_G$ can be composed with $\Gamma_G\colon (0,\infty]^{\vec E}\to[0,\infty)^{\vec E}$, yielding a two-step local update rule on $(0,\infty]^{\vec E}$ which will now play a crucial role. 
\begin{theorem}[The cavity equation at infinite activity]
\label{th:infiniteactivity}
Let $G$ be a cavity-monotone network on which the cavity equation at activity $t$ admits a unique globally attractive fixed point $\bx(t)$, for every $t>0$. Then,
\begin{equation}
\label{eq:uxox}
\overline{\m}:=\lim_{t\to\infty}\uparrow \m(t)
\end{equation}
exists in $(0,\infty]^{\vec E}$, and is the smallest solution to the so-called \textit{cavity equation at infinite activity} on $G$ :
\begin{equation}
\label{eq:reczero}
\overline{\m}=\left(\overline{\Gamma}_G\circ\Gamma_G\right)(\overline{\m}).
\end{equation}
In particular, for any unimodular probability measure $\cL$ concentrated on cavity-monotone networks, and satisyfing $\cL\left[\rm{rank}(\mu_\circ)\right]<\infty$, we have 
$$m(\cL)=\frac{1}{2}\cL\left[U_{\mu_\circ}\left(\xo_{i\to\circ}:i\in\partial\circ\right)\right].$$
\end{theorem}
\bp{ of Theorem \ref{th:infiniteactivity}}
By assumption $\bx^{k}(t)\to\bx(t)$ for any $t>0$, where $\bx^0\equiv\bf{0}$ and for all $k\in\N$,
\begin{eqnarray}
\label{eq:klambda}
\bx^{k+1}(t)=t\Gamma_G\left(t\frac{\bx^k(t)}{t}\right)
& \textrm{ or equivalently, } & 
\frac{\bx^{k+1}(t)}{t}=\Gamma_G\left(\bx^{k}(t)\right) .
\end{eqnarray}
But $(t,\bx)\mapsto t\Gamma_G(t\bx)$ is increasing in $t$ and decreasing in $\bx$, so an immediate induction over $k$ shows that $t\mapsto t^{-1}{\bx^{k}(t)}$ and $t\mapsto \bx^{k}(t)$ are respectively non-increasing and non-decreasing, for every $k\in\N$. Thus, $t\mapsto \bx(t)$ is non-decreasing, hence the existence of (\ref{eq:uxox}). The identity (\ref{eq:reczero}) is then obtained by passing to the limits in (\ref{eq:klambda}). Finally, if $\by\in(0,\infty]^{\vec{E}}$ satisfies $\by=(\overline{\Gamma}_G\circ\Gamma_G)(\by)$, then for every $k\in\N$ and $t>0$,
\begin{equation}
\label{eq:maximal}
\bx^{2k}(t) \leq \by.
\end{equation}
Indeed, the property is trivial when $k=0$, and is preserved from $k$ to $k+1$ because $t\Gamma_G(t\bw)\leq\overline{\Gamma}_G(\bw)$ for any $\bw\in [0,\infty)^{\vec{E}}$.
Letting $k\to\infty$ and then $t\to\infty$ in (\ref{eq:maximal}) yields $\bxo\leq\by$, as desired. 
\ep
\begin{remark}[Continuity with respect to the activity]
\label{re:continuity}
Incidentally, we have just obtained that $t\mapsto t^{-1}{\bx(t)}$ and $t\mapsto \bx(t)$ are respectively non-increasing and non-decreasing, so that 
\begin{eqnarray*}
0< s\leq t & \Longrightarrow & \frac{s}{t}\bx(t) \leq \bx(s)\leq \bx(t).
\end{eqnarray*}
This guarantees the continuity of $t\mapsto \bx(t)$, and hence that of $u(\cL;t)$, as promised in Remark \ref{re:ldp}. 
\end{remark}

\subsection{Explicit resolution for $b-$matchings on GW trees}

Many classical sequences of diluted random graphs 
%(including Erd\H{o}s-R\'enyi graphs with fixed average degree and random graphs with a prescribed degree profile) 
admit almost surely a particularly simple random weak limit $\cL$, namely a \textit{unimodular Galton-Watson (UGW) tree} (see Example 1.1 in \cite{aldlyo}). This random rooted tree is parametrized by a probability distribution $\pi\in\cP(\N)$ with finite mean, called its \textit{degree distribution}. It is obtained by a Galton-Watson branching process where the root has offspring distribution $\pi$ and all other genitors have the  size-biased offspring distribution $\widehat{\pi}\in\cP(\N)$ defined by 
\begin{equation*}
\label{eq:F}
\forall n\in \N, \widehat{\pi}_{n} = {(n+1) \pi_{n+1}}/{\sum_{k} k \pi_k}.
\end{equation*}
Thanks to the markovian nature of the branching process, the cavity equation at infinite activity simplifies into a \textit{recursive distributional equation} (RDE) (see \cite{aldousbandyopadhya} for a survey). Let us describe it and solve it in the case of $b-$matchings, where $b\geq 1$ is fixed. The local cavity and energy ratios are simply
\begin{eqnarray*}
\Gamma(x_1,\ldots,x_n)=\frac{\sum_{I\subseteq [N]:|I|\leq b-1}\prod_{i\in I}x_i}{\sum_{I\subseteq [N]:|I|\leq b}\prod_{i\in I}x_i}
& \textrm{ and } &
U(x_1,\ldots,x_n) = \frac{\sum_{I\subseteq [N]:|I|\leq b}|I|\prod_{i\in I}x_i}{\sum_{I\subseteq [N]:|I|\leq b}\prod_{i\in I}x_i}.
\end{eqnarray*}
In the infinite activity limit, the local cavity ratio becomes 
\begin{eqnarray*}
\overline{\Gamma}(x_1,\ldots,x_n) & = & \frac{\sum_{I\subseteq [N]:|I|= b-1}\prod_{i\in I}x_i}{\sum_{I\subseteq [N]:|I|= b}\prod_{i\in I}x_i},
\end{eqnarray*}
where all conventions regarding degenerate cases are obtained by taking the corresponding limits. Given $Q\in\cP((0,\infty])$, we let $\Theta(Q)\in\cP\left([0,1])\right)$ denote the law of  $\Gamma\left(Y_{1},\ldots,Y_{\widehat{N}}\right)$, where $\widehat{N}$ has law $\widehat{\pi}$ and $Y_1,Y_2,\ldots$ are i.i.d. with law $Q$,  independent of $\widehat N$. Similarly, given $P\in\cP([0,1])$, we let 
$\overline{\Theta}(P)\in\cP\left((0,\infty]\right)$ denote the law of  $\overline{\Gamma}\left(X_{1},\ldots,X_{\widehat{N}}\right)$, where $\widehat{N}$ has law $\widehat{\pi}$ and $X_1,X_2,\ldots$ are i.i.d. with law $P$,  independent of $\widehat N$. 
Thanks to the markovian nature of the GW branching process, the law $Q\in\cP((0,\infty])$ of a message sent towards the root in the configuration $\bxo$ must satisfy the RDE $Q  = (\overline{\Theta}\circ\Theta)(Q)$. Equivalently, $P=\Theta(Q)$ must satisfy
$P  = (\Theta\circ\overline{\Theta})(P).$ More precisely, letting $M(P)$ denote the expectation of $\frac{1}{2}U\left(Y_{1},\ldots,Y_{N}\right)$, where $N$ has law $\pi$ and $Y_1,Y_2,\ldots$ are i.i.d. with law $Q=\overline{\Theta}(P)$, independent of $N$, we may rephrase Theorem \ref{th:infiniteactivity} as follows (distributions are here endowed with the usual stochastic order) :
\begin{corollary}
$m(\cL)=M(P)$, where $P\in\cP([0,1])$ is the smallest solution to the RDE $P  = (\Theta\circ\overline{\Theta})(P)$.
\end{corollary}
The fixed points of $\Theta\circ\overline{\Theta}$ turn out to be in one-to-one correspondance with the \textit{historical minima} of a certain function $H\colon[0,1]\to\R$ defined in terms of the degree generating function $\phi(s)=\sum_{k}\pi_ks^k$ :
\begin{equation*}
H(s)=b-\frac{b}{2}g(s)-\frac{b}{2}(g\circ f)(s)+\frac{c}{2}f(s)(f\circ f)(s),
\end{equation*}
\begin{equation*}
\textrm{ with }\qquad c=\phi'(1),\qquad\qquad f(s)=\frac{1}{c}\sum_{k=0}^{b-1}\frac{s^{k}\phi^{(k+1)}(1-s)}{k!} \qquad\textrm{ and }\qquad  g(s)=\sum_{k=0}^b\frac{t^{k}\phi^{(k)}(1-s)}{k!}.
\end{equation*}
A \textit{historical minima} of $H$ is a number $s\in [0,1]$ satisfying $H'(s)=0$ and $H(t)>H(s)$ for all $t\in[0,s)$. 
%Note that there are only finitely many such numbers, since $H$ is analytic and $[0,1]$ is compact.
\begin{theorem}[Resolution of the RDE]
\label{th:rde}
Let $s_1<\ldots<s_r$ denote the historical minima of the function $H$. Then, the distributional equation $P  = (\Theta\circ\overline{\Theta})(P)$ admits exactly $r$ solutions, and they are stochastically ordered : $P_1<\ldots<P_r$. Moreover, for each $1\leq i\leq r$, we have $M(P_i)=H(s_i)$.  
\end{theorem}
In particular, $m(\cL)=\min H$, which is exactly the formula given in Theorem \ref{th:KSb}. Theorem \ref{th:rde} was established in \cite{rankarxiv} for the special case $b=1$, but the proof can easily be adapted to the general case. For the sake of completeness, we have included a general proof in the Appendix. 

\section*{Appendix : resolution of the RDE}
First observe that the mappings $\Theta, \overline{\Theta}$ and $M$ are all decreasing with respect to stochastic order, and continuous with respect to the topology of weak convergence. Note also that $c\widehat{f}'(t)t=bf(t)$, so that 
$H'(t)=c f'(t)\left((f\circ f)(t)-t\right).$ Thus, $H'(t)=0$ if and only if $(f\circ f)(t)=t.$

\begin{lemma}[Properties of $\Gamma,\overline{\Gamma}$ and $U$]\mbox{}
\label{lm:propgamma}
\begin{enumerate}
\item Let $(x_1,\ldots,x_n) \in [0,1]^n$. Then,
\begin{enumerate}
\item $\overline{\Gamma}(x_1,\ldots,x_n)=\infty \Longleftrightarrow \sum_{i=1}^{n}\ind_{\left\{x_i>0\right\}} < b$ ;
\item Setting $y_i=\overline{\Gamma}(x_k:k\neq i)$, we have $\displaystyle{\sum_{i=1}^n\frac{x_iy_i}{1+x_iy_i}\ind_{\{y_i<\infty\}} = b\ind_{\left\{\sum_{i=1}^{n}\ind_{\left\{x_i>0\right\}}>b\right\}}}$
\end{enumerate}
\item Let $(y_1,\ldots,y_n) \in (0,\infty]^n$. Then,
\begin{enumerate}
\item $\Gamma(y_1,\ldots,y_n)>0 \Longleftrightarrow \sum_{i=1}^{n}\ind_{\left\{y_i=\infty\right\}}<b$;
\item Setting $x_i'=\Gamma(y_k:k\neq i)$, we have $\displaystyle{\sum_{i=1}^n\frac{x_i'y_i}{1+x_i'y_i}\ind_{\{x_i'<\infty\}} = U(x_1',\ldots,x_n') - b\wedge \sum_{i=1}^{n}\ind_{\left\{y_i=\infty\right\}}}.$
%.\\
\end{enumerate}
\end{enumerate}
\end{lemma}
\bp{}
Properties $1.a$ and $2.a$ are straightforward from the definition of $\overline{\Gamma}$ and $\Gamma$. Regarding property $1.b$, set $K=\#\left\{i\in[n]:x_i>0\right\}.$ If the sum is non-zero then there must be an $i$ satisfying both $y_i>0$ and $x_i<\infty$. By $1.a$, this implies $K>b$. Conversely, if $K>b$ then $x_i<\infty$ for every $i\in[n]$. We have just shown
\begin{eqnarray*}
\sum_{i=1}^n\frac{x_iy_i}{1+x_iy_i}\ind_{\{y_i<\infty\}} & = & \ind_{\left\{K>b\right\}}\sum_{i=1}^n\frac{x_iy_i}{1+x_iy_i}.\\
& = & \ind_{\left\{K>b\right\}}\sum_{i=1}^n\frac{\sum_{|I|=b,I\ni i}\prod_{k\in I}x_k}{\sum_{|I|=b}\prod_{k\in I}x_k}\\
& = & b\ind_{\left\{K> b\right\}},
\end{eqnarray*}
where the second equality is obtained by replacing $y_i=\overline{\Gamma}(x_k:k\neq i)$ by its definition. 
For property $2.b$, set $L=\#\left\{i\in[n]:y_i=\infty\right\}.$ When $L=0$, $2.b$ boils down to formula (\ref{eq:energycavity}). The case $1\leq L \leq b$ can then be obtained by successively setting each of the $L$ variables to $\infty$, which amounts to condition on the presence of the corresponding ground elements. For $L\geq b$, both sides of the equation are zero. 
\ep
\begin{lemma}
\label{lm:last}
Assume that $P\xrightarrow{\overline{\theta}}Q\xrightarrow{\theta}P'$. Set $s=P(\{0\}^c)$, $t=Q(\{\infty\})$ and $s'=P'(\{0\}^c)$. Then,
\begin{enumerate}
\item $s \xrightarrow{f} t \xrightarrow{f} s'$ ;
\item $P'\leq P \Longrightarrow M(P)\leq H(s)$ ;
\item $P'\geq P \Longrightarrow M(P)\geq H(s)$.
\end{enumerate}
In particular, if $P  = (\Theta\circ\overline{\Theta})(P)$ then $M(P) = H(s)$ and $H'(s)=0$.
\end{lemma}
\bp{}
In the whole proof, $N$ denotes a generic random integer with law  $\pi$, $\widehat{N}$ a generic random integer with law $\widehat{\pi}$, $X,X_1,X_2,\ldots$ generic $(0,\infty]-$valued random variables with law $P$, $Y,Y_1,Y_2,\ldots$ generic $[0,1]-$valued random variables with law $Q$, and $X',X_1',X_2',\ldots$ generic $(0,\infty]-$valued random variables with law $P'$. 
We use the convention that all variables appearing under the same expectation are independent. 
With these notations, properties $1.a$ and $2.a$ in Lemma \ref{lm:propgamma} give 
$$\PP\left(Y=\infty\right) = \PP\left(\sum_{i=1}^{\widehat N}\ind_{\left\{X_i>0\right\}} < b\right)\quad\textrm{ and }\qquad\PP\left(X'>0\right) = \PP\left(\sum_{i=1}^{\widehat N}\ind_{\left\{Y_i=\infty\right\}}< b\right),$$
which, in view of the definition of $f$, yields exactly the first claim of the Lemma. Now, using property $1.b$ and $2.b$, we respectively obtain the two following identities :
\begin{eqnarray*}
\EE\left[\frac{XY}{1+XY}\ind_{\{Y<\infty\}}\right] & = &
\sum_{n\in\N}\widehat{\pi}(n)\EE\left[\frac{X\overline{\Gamma}\left(X_1,\ldots,X_{n-1}\right)}{1+X\overline{\Gamma}\left(X_1,\ldots,X_{n-1}\right)}\ind_{\left\{\overline{\Gamma}\left(X_1,\ldots,X_{n-1}\right)<\infty\right\}}\right]\\
& = &
\sum_{n\in\N}\pi(n)n\EE\left[\frac{X\overline{\Gamma}\left(X_1,\ldots,X_{n-1}\right)}{1+X\overline{\Gamma}\left(X_1,\ldots,X_{n-1}\right)}\ind_{\left\{\overline{\Gamma}\left(X_1,\ldots,X_{n-1}\right)<\infty\right\}}\right]\\
& = & \sum_{n\in\N}\pi(n)\EE\left[\sum_{i=1}^n\frac{X_i\overline{\Gamma}\left(X_k:k\neq i\right)}{1+X_i\overline{\Gamma}\left(X_k:k\neq i\right)}\ind_{\left\{\overline{\Gamma}\left(X_k : k\neq i\right)<\infty\right\}}\right]\\
 & = & b\PP\left(\sum_{i=1}^N\ind_{\left\{X_i>0\right\}}> b\right)\\
 & = & b(1-g(s)).\\
 \\
%\textrm{Similarly, using property 2.b gives}\\
\\
\EE\left[\frac{X'Y}{1+X'Y}\ind_{\{Y<\infty\}}\right] & = & \sum_{n\in\N}\widehat{\pi}(n)\EE\left[\frac{Y\Gamma\left(Y_1,\ldots,Y_{n-1}\right)}{1+Y\Gamma\left(Y_1,\ldots,Y_{n-1}\right)}\ind_{\{Y<\infty\}}\right]\\
& = &  \sum_{n\in\N}\pi(n)n\EE\left[\frac{Y\Gamma\left(Y_1,\ldots,Y_{n-1}\right)}{1+Y\Gamma\left(Y_1,\ldots,Y_{n-1}\right)}\ind_{\{Y<\infty\}}\right]\\
& = & \sum_{n\in\N}\pi(n)\EE\left[\sum_{i=1}^n\frac{Y_i\Gamma\left(Y_k:k\neq i\right)}{1+Y_i\Gamma\left(Y_k:k\neq i\right)}\ind_{\{Y_i<\infty\}}\right]\\
& = & \EE\left[U(Y_1,\ldots,Y_N)\right] - \EE\left[b\wedge\sum_{i=1}^N\ind_{\left\{Y_i=\infty\right\}}\right] \\
& = & 2M(P) - b(1-g(s)) - csf(s).\\
\end{eqnarray*}
Since the mapping $\frac{xy}{1+xy}$ is increasing in $x$, claims 2 and 3 follow.
\ep

\bp{ of  Theorem \ref{th:rde}} Fix $s\in[0,1]$ satisfying $H'(s)=0$, i.e. $(f\circ f)(s)=s$. Define $P_0^{s}\in\cP([0,1])$ to be the Bernoulli distribution with parameter $s$, and set then iteratively $$P_{n+1}^{s}=(\Theta\circ\overline{\Theta})(P_{n}^{s})$$ for all $n\in\N$. By part 1 in Lemma \ref{lm:last}, $P_1^{s}$ is a distribution on $[0,1]$ satisfying $P_1(\{0\}^c)=s$. Since $P^{s}_0$ is the largest such distribution, we have $P^{s}_1\leq P^{s}_0$. But both $\Theta$ and $\overline{\Theta}$ are decreasing, so by immediate induction, the sequence $(P_n^{s})_{n\in\N}$ is non-increasing. Thus, the limit $P^{s}=\lim_{n\to\infty}\downarrow P_n^{s}$ exists in $\cP([0,1])$. Setting $s_\infty=P^{s}(\{0\}^c)$, we claim that
\begin{enumerate}
\item $P^{s}$ is a fixed point of $\Theta\circ\overline{\Theta}$ ;
\item $M(P^{s})=H(s_\infty);$
\item $s_\infty\leq s;$
\item $H(s_\infty)\leq H(s)$.
\end{enumerate}
Part 1 follows from the continuity of $\Theta$ and $\overline{\Theta}$. Part 2 is guaranteed by Lemma \ref{lm:last}. Part 3 is a consequence of the fact that $P^{s} \leq P^{s}_0$. Finally, for each $n\in\N$, we have $P_n(\{0\}^c)=s$ and $P^{s}_n\geq P^{s}_{n+1}$, so Lemma \ref{lm:last} guarantees that $M(P^{s}_n) \leq H(s)$. Letting $n\to\infty$ yields exactly part 4. 

We are now in position to prove the equivalence between the historical minima of $H$ and the solutions to $(\Theta\circ\overline{\Theta})(P)=P$. If $s$ is a historical minimum, then parts 3 and 4 force $s_\infty = s$ so $P=P^{s}$ is a solution satisfying $M(P)=H(s)$, as desired. Conversely, we must show that any solution $P$ is in fact of the form $P^s$ for some historical minimum $s$. Set $s=P(\{0\}^c)$, which satisfies $H'(s)=0$ by Lemma \ref{lm:last}. We have $P\leq P^s_n$ for any $n\in\N$ : this holds for $n=0$ because $P^s_0$ is the largest element of $\cP([0,1])$ such that $P_0(\{0\}^c)=s$, and it then inductively extends to all $n\in\N$ by monotony of $\Theta\circ\overline{\Theta}$. Letting $n\to\infty$, we obtain $P\leq P^s$ ; but by Lemma \ref{lm:last} we also have $M(P)=H(s)=M(P^s)$. Thus, $P=P^s$ ($M$ is decreasing). Finally, if $t<s$ is any historical minimum then clearly $P^t_0\leq P^s_0$, which implies $P^t\leq P^s$. In fact the inequality is strict, because $P^t(\{0\}^c)=t<s=P^s(\{0\}^c)$. Applying the decreasing mapping $M$ yields $H(t)>H(s)$, which shows that $s$ is a historical minimum.

\bibliographystyle{abbrv}
\bibliography{cavity}
\end{document}